# DISCUSSION: THE DANTZIG SELECTOR: STATISTICAL ESTIMATION WHEN $p$ IS MUCH LARGER THAN $n$[1]

By Peter J. Bickel

*University of California, Berkeley*

**1. A personal comparative review.** Reading this very interesting paper prompted me to review the current literature on the Lasso and sparsity and realize that there are at least three different sets of conditions on the behavior of the predictors $\mathbf{X}_1, \ldots, \mathbf{X}_p$ under which the coefficients of the sparsest representation of a general regression model are recovered with an $l_2$ error of order $\sqrt{\frac{s}{n} \log p}$, where $s$ is the dimension of the sparsest model. These are, respectively, the conditions of this paper using the Dantzig selector and those of Bunea, Tsybakov and Wegkamp [2] and Meinshausen and Yu [9] using the Lasso. Strictly speaking, Bunea, Tsybakov and Wegkamp consider only prediction, not $l_2$ loss, but in a paper in preparation with Ritov and Tsybakov we show that the spirit of their conditions is applicable for $l_2$ loss as well. Since these authors emphasize different points and use different normalizations, I thought it would be useful to present them together. Write the model as

$$Y = X_{n \times p} \boldsymbol{\beta} + \boldsymbol{\varepsilon},$$

where for simplicity we take $\boldsymbol{\varepsilon} \sim N(0, \sigma^2 I_n)$, a case falling under all the authors' conditions, and

$$X = (\mathbf{X}_1, \ldots, \mathbf{X}_p).$$

We begin by assuming that

$$X\boldsymbol{\beta} \in [\mathbf{X}_{i_1}, \ldots, \mathbf{X}_{i_s}] \equiv V,$$

an unknown unique $s$-dimensional linear subspace of $[\mathbf{X}_1, \ldots, \mathbf{X}_p]$ and no lower-dimensional subspace. For simplicity, we follow Knight and Fu [8] and Meinshausen and Yu [9] and assume $|\mathbf{X}_j|^2 = n$, $1 \leq j \leq p$, which puts the problem on the same scale as the familiar: $X$ a matrix of $n$ i.i.d. $p$ vectors.

Received February 2007.
[1]Supported in part by NSF Grant DMS-06-05236.

In its current form, a goal of all three authors is to state conditions on $X$ under which, for the method of estimation they propose,

$$(1) \qquad |\hat{\boldsymbol{\beta}} - \boldsymbol{\beta}^\circ| = O_p\left(\sqrt{\frac{s}{n}\log p}\right),$$

where $\boldsymbol{\beta}^\circ$ corresponds to the "true" $\boldsymbol{\beta}$. This can be done by thresholding if $\mathbf{X}_1, \ldots, \mathbf{X}_p$ are orthogonal and is optimal; see the work of Donoho and Johnstone [5].

In general, thresholding the least squares estimates of $\boldsymbol{\beta}$ is incorrect since when $X^T X$ has rank $n < p$, the LSE are not uniquely defined and the Dantzig selector and Lasso are computationally feasible ways of achieving (1) under different sets of assumptions. In the orthogonal case they both correspond to soft thresholding, as Candès and Tao [3] point out. All of these are based on the two key geometric qualities introduced by Donoho and Candès and Tao.

Following Donoho [4] and Candès and Tao [3], define $\varphi_{\min}(m)$ to be the minimal eigenvalue of the Gram matrix $X_L^T X_L$ for all $L$ with $|L| \leq m$ and $\varphi_{\max}$ to be the maximum eigenvalue of $X^T X$. Here, $X_L = \{(X_{i_1}, \ldots, X_{i_L}) : i_j \in L\}$. Let

$$\theta_{m,m'} = \max\{|(X_L \mathbf{c}_L, X_L \mathbf{c}_{L'})| : L \cap L' = \varnothing,$$
$$|L| \leq m', |L'| \leq m, |\mathbf{c}_L| \leq 1, |\mathbf{c}_{L'}| \leq 1\}$$

where $\mathbf{c}_L$ is $|L| \times 1$ and $\mathbf{c}_{L'}$ is $|L'| \times 1$.

Denote the Lasso estimate by

$$\hat{\boldsymbol{\beta}}_L \equiv \arg\min\left\{|\mathbf{Y} - X\boldsymbol{\beta}|^2 + \lambda \sum_{j=1}^p |\beta_j|\right\}.$$

Then the conditions common to both the Dantzig selector and the Lasso are:

A1. $\varphi_{\max} \leq \overline{k} < \infty$.

The unicity of sparsest representation condition:

A2. $\varphi_{\min}(2s) \geq \underline{k} > 0$.

The differences come in the condition which goes beyond A2 specifying how "near orthogonality" $S$ and $\bar{S}$ are, where $\bar{S}$ is the complement of $S$.

The condition of Bunea, Tsybakov and Wegkamp [2] is:

A3 (BTW). $\rho_s \equiv \max\{|(\mathbf{X}_i, \mathbf{X}_j)| : i \in L, j \in \bar{L}, |L| \leq s\} \leq \frac{M}{s}$ for a constant $M \leq \frac{1}{32}$.

That of Meinshausen and Yu [9] is a great strengthening of A2:

A3 (MY). $\varphi_{\min}(s \log n) \geq \varepsilon > 0$ for some $\varepsilon > 0$.

The Candès and Tao [3] condition is, in these terms:

A3 (CT). $\theta_{s,2s} < \varphi_{\min}(2s) < 1$ and $\varphi_{\max}(2s) + \theta_{s,2s} < 2$.

That $\rho_s$ must decrease at rate $\frac{d}{s}$ is clearly meaningful. As $V$ gets bigger, all columns not in $V$ need to become more and more orthogonal to it. How are these conditions related? I hope that the other discussants will shed light on this.

**2. The form of the Dantzig predictor.** There is a compelling reason for using $\|X^T(\mathbf{Y} - X\boldsymbol{\beta})\|_\infty$ rather than $\|\mathbf{Y} - \mathbf{X}\beta\|_\infty$ as part of the objective, which the authors are undoubtedly aware of, hinted at in Section 1.3, but did not make explicit. Suppose the $(\mathbf{X}^i, \mathbf{Y}_i)$, $1 \leq i \leq n$, are i.i.d. where $\mathbf{X}^i$ is $p$-dimensional and we want to estimate $f(\mathbf{X}) = E(Y|\mathbf{X})$ using a dictionary $\{f_j\}$, $j \geq 1$. Usually we imagine that $f(\mathbf{X}) = \sum_{j=1}^\infty \beta_j f_j(X)$ where the expansion is in the $L_2$ sense and we expect that $\sum_{j=1}^p \beta_j f_j(\mathbf{X})$ is a good approximation in the $L_2$ sense. If we now identify $\mathbf{X}_j$ with $(f_j(\mathbf{X}_1), \ldots, f_j(\mathbf{X}_n))$, then the minimizer of $\|\mathbf{Y} - X\boldsymbol{\beta}\|_\infty$ is trying to minimize $\operatorname{ess\,sup} |f(\mathbf{x}) - \sum_{j=1}^p \beta_j f_j(\mathbf{x})|$ which can, no matter what $p$ is, be very large. On the other hand, the Dantzig predictor needs to match correctly all $\int f(x) f_j(x) \, dP(x)$, $1 \leq j \leq p$, which is just what is wanted.

**3. Model selection.** As Candès and Tao [3] point out, from a statistical point of view the Dantzig selector, just as the Lasso, can be viewed as a method of model selection. Of course, the $l_2$ norm results per se do not say which variables are the ones appearing in the sparsest model. But as they, and more explicitly, Meinshausen and Yu [9] point out, model selection after computation of the estimate (Dantzig selector or Lasso, resp.) can give an idea of what are the variables with large true coefficients; see also Wainwright [10] and Zhao and Yu [11]. I would argue that, in the large $p$, large $n$ situations being considered, there may well be variables which are interpretable and, though necessarily highly correlated with variables which do appear in the sparsest model, themselves do not appear. Or, even worse, condition A2 may not hold. There may be two spaces $S$ and $S'$ of the same dimension $s$, each of which provides sparsest representations. Both the Dantzig selector and the Lasso in a case like this will, depending on starting point, converge to a point close to one of the representations—but in fact we would like to consider the whole space of such representations.

Here is a simple example with $s = 3$. For simplicity, we consider the noiseless case. For simplicity, make the rows of $X$ i.i.d. so that we can talk in population terms. We assume $X_1$, $X_2$, $X_3$ are predictors with mean 0 and variance 1. $X_1$ and $X_3$ are independent and

$$X_2 = \alpha X_1 + \beta X_3.$$

Suppose, also, $Y = X_1 + X_2 + X_3$. Then, there are three representations of $Y$ corresponding to $(X_1, X_2)$, $(X_2, X_3)$ and $(X_1, X_3)$, and a Lasso fit might converge to any of the three. If the number of predictors is small, exact collinearities or near collinearities can be identified, but this is not the case when $p$ is large.

Suppose though that we can identify collinearities as above. Is there a reasonable measure of importance of a variable? In the above example there are two representations involving each variable. For instance, for $X_2$,

$$Y = \left(1 + \frac{1}{\beta}\right)X_2 + \left(1 - \frac{\alpha}{\beta}\right)X_1,$$

$$Y = \left(1 + \frac{1}{\alpha}\right)X_2 + \left(1 - \frac{\beta}{\alpha}\right)X_3.$$

Rewrite the models so as to take out the effect of the other variable on $X_2$:

$$Y = \left(1 + \frac{1}{\beta}\right)(X_2 - \alpha X_1) + (\alpha + 1)X_1,$$

$$Y = \left(1 + \frac{1}{\alpha}\right)(X_2 - \beta X_3) + (\beta + 1)X_1.$$

A natural measure of the importance of $X_2$ in the presense of $X_1$ is then

$$SN^2(2|1) \equiv \left(1 + \frac{1}{\beta}\right)^2 \beta^2 = (\beta + 1)^2,$$

and similarly,

$$SN^2(2|3) \equiv (\alpha + 1)^2.$$

It seems reasonable that we ascribe to $X_2$ overall importance as

$$I(2) = \max\{(\alpha + 1)^2, (\beta + 1)^2\}.$$

If we add independent noise $\varepsilon_i$ with variance $\sigma^2$ to each observation,

$$Y = X_{1i} + X_{2i} + X_{3i} + \varepsilon_i,$$

then the natural definition of importance is

$$I(2) = \max\left\{\frac{SN^2(2|1)}{\sigma^2}, \frac{SN^2(2|3)}{\sigma^2}\right\}.$$

If we estimate the coefficients of $X_2$ in each of the two models, we see that up to a constant, $\hat{I}(2)$ is just the maximum of the squares of the $t$ statistics for testing the hypothesis that the coefficient of $X_2$ is 0 in each of the two models.

How should we address the question of determining "importance" of a predictor? What we need first is the linear space spanned by $\{\mathbf{X}_{i_j} : j \in V\}$ of all vectors $\mathbf{X}_j$ uncorrelated with $X\boldsymbol{\beta}$. This can be done easily by using the test statistic $\frac{\mathbf{X}_k^T \mathbf{Y}}{\sqrt{n}\hat{\sigma}}$ and using a cut-off of $\sqrt{\frac{2\log p}{n}}$ where $\hat{\sigma}^2$ is the residual mean square after a best predictive Lasso fit. This approach is worked out to eliminate more than just variables uncorrelated with $X\boldsymbol{\beta}$ by Fan and Lv [6], who focus on situations where $p \gg n$. Once we have $V$, we would like to assess the importance of our variables. If $|V| = s$, we are done using the usual least squares fit. However, we would typically have $|V| \gg s$. In that case since every variable can enter into a representation $X\boldsymbol{\beta}$, we need to apply the Dantzig selector or the Lasso to reduce dimension. Unfortunately, unless all representations have $|\boldsymbol{\beta}|_1$ which differs by at most $o_p(\sqrt{\log p/n})$, both the Dantzig selector and the Lasso will converge to a minimum $|\boldsymbol{\beta}|_1$ solution. It is not clear to me how to obtain alternative representations when there are many variables without considering all possible subspaces, which is hopeless. However, it seems plausible that a first approximation to "importance" might be obtained as follows. Select a set of $K$ candidate variables, say, the $K$ most correlated with $Y$ whether they appear or not in the Lasso or Dantzig fit. Regress each of these together with all but one of the $N$ variables having nonzero coefficients in the Dantzig or Lasso fit on the fitted values obtained by Dantzig or Lasso. Retain only candidate variables such that, in at least one of the $N$ resulting LS fits, the fit is good and the importance of the candidate in that context is large.

The properties of such procedures under the circumstances I outline, many variables and collinearity, seem worth considering. Of course, not all nonzero coefficient variables in the Dantzig or Lasso fits necessarily have to be considered. Again one could limit to a candidate set having large coefficients.

**4. Choice of $\lambda$ and $\sigma$.** Candès and Tao [3] do not dwell on the choice of $\lambda$, which should be of order $\sigma\sqrt{2\log n/p}$, but, in Section 4, simply specify the empirical maximum over several realizations of $|X^T \mathbf{Z}|$, or in our case $\frac{|X^T \mathbf{Z}|}{\sqrt{n}}$ with $\mathbf{Z} \sim \mathcal{N}(\mathbf{0}, I_n)$. If the goal is to minimize the $l_2$ norm of $|\hat{\boldsymbol{\beta}}_{\mathrm{CT}} - \boldsymbol{\beta}_0|$, where $\hat{\boldsymbol{\beta}}_{\mathrm{CT}}$ is their estimate, it is worth pointing out that $V$-fold cross-validation can be used here as well. That is, choose a test subsample of size, say, $\log n$ and fit the Dantzig predictor to the remaining $n - \log n$ observations, say, for a grid of $\lambda$ values of width $c\sqrt{\frac{\log n}{p}}$. Then choose the best $\lambda$ in the sense of $l_2$ prediction for the test sample. This should give optimal rates if $\log n = o(\log p)$, the typical case. If all elements are bounded in absolute value, this follows from Theorem 6 of Bickel, Ritov and Zakai [1]; see also Györfi, Kohler, Krzyżak and Walk [7]. Strictly speaking, this optimizes for prediction loss, not $l_2$, but for $\boldsymbol{\beta}_S$ the two are equivalent.

## REFERENCES

[1] BICKEL, P., RITOV, Y. and ZAKAI, A. (2006). Some theory for generalized boosting algorithms. *J. Mach. Learn. Res.* **7** 705–732. MR2274384

6 P. J. BICKEL

[2] Bunea, F., Tsybakov, A. B. and Wegkamp, M. H. (2007). Aggregation for Gaussian regression. *Ann. Statist.* **35** 1674–1697. MR2351101

[3] Candès, E. J. and Tao, T. (2007). The Dantzig selector: Statistical estimation when $p$ is much larger than $n$. *Ann. Statist.* **35** 2313–2351.

[4] Donoho, D. L. (2006). For most large underdetermined systems of linear equation the minimal $l_1$-norm solution is also the sparsest solution. *Comm. Pure Appl. Math.* **59** 797–829. MR2217606

[5] Donoho, D. L. and Johnstone, I. M. (1994). Ideal spatial adaptation by wavelet shrinkage. *Biometrika* **81** 425–455. MR1311089

[6] Fan, J. and Lv, J. (2006). Sure independence screening for ultra-high dimensional feature space. Unpublished manuscript. Available at arxiv.org/abs/math/0612857.

[7] Györfi, L., Kohler, M., Krzyżak, A. and Walk, H. (2002). *A Distribution-Free Theory of Nonparametric Regression*. Springer, New York. MR1920390

[8] Knight, K. and Fu, W. (2000). Asymptotics for lasso-type estimators. *Ann. Statist.* **28** 1356–1378. MR1805787

[9] Meinshausen, N. and Yu, B. (2007). Lasso-type recovery of sparse representations for high-dimensional data. *Ann. Statist.* To appear.

[10] Wainwright, M. (2006). Sharp thresholds for high-dimensional and noisy sparsity recovery using $l_1$-constrained quadratic programs. In *Proc. Allerton Conference on Communication, Control and Computing.*

[11] Zhao, P. and Yu, B. (2007). Model selection with the Lasso. *J. Mach. Learn. Res.* To appear.

Department of Statistics
University of California, Berkeley
Berkeley, California 94720-3860
USA
E-mail: bickel@stat.berkeley.edu